# Non-homogeneous Linear Set-Valued Differential Equations with Variable Matrix Coefficients

Uma Maheswara Rao Epuganti[1] and Gnana Bhaskar Tenali[2]

Department of Applied Mathematics and Systems Engineering

Florida Institute of Technology

Melbourne, FL 32901, USA

November 20, 2025

**Abstract**

We investigate the initial value problems for non-homogeneous linear differential equations whose solutions are set-valued maps taking values in the space of nonempty compact convex subsets of $\mathbb{R}^2$, denoted by $K_c(\mathbb{R}^2)$. The differential formulation is based on the generalized derivative that includes the Hukuhara derivative, as well as its extensions, Bede-Gal (BG), and Plotnikov-Skripnik (PS) derivatives, and we obtain some general as well as constructive formulas for the solutions. Several illustrative examples are provided.

## 1 Introduction

The study of differential equations for (compact and convex) set-valued maps was initiated in [1], using the notion of Hukuhara derivative [2]. The monograph [3] summarizes the evolution of ideas in the study of set-valued differential equations (SVDEs) in different directions. Over the past several years, attempts have been made to rectify the intrinsic disadvantage of using the Hukuhara derivative in the study of SVDEs, namely that any solution of such a SVDE has a nondecreasing diameter. In an effort to fix this issue in the context of the study of fuzzy differential equations, [4, 5] introduced the notion of generalized differentiability, that is now known as Bede-Gal differentiability (BG-differentiability), and [6] introduced the notion of Plotnikov-Skripnik differentiability (PS-differentiability). Much work has been done to study the properties of these derivatives (see [7–12]). For an overview of the evolution of some of these ideas in the context of fuzzy differential equations and differential inclusions, we refer the reader to [13] and [14].

In recent years, significant attention has also been directed towards the qualitative analysis of SVDEs with the Hukuhara derivative, towards introducing novel techniques for analyzing stability, constructing Lyapunov functions, and controlling set-valued dynamical systems using geometric methods. Notable contributions include the generalization of the Liouville formula, impulsive stabilization of ensembles, and stability analysis in both continuous and discrete settings (see [15–17]).

Several fundamental results concerning linear homogeneous SVDEs are studied in a series of papers [7, 9, 10, 12]. In order to study such SVDEs from a control theory perspective, there is a need for a systematic study of non-homogeneous SVDEs under the notion of generalized differentiability. Following the approach taken in [10], we studied non-homogeneous linear SVDEs with matrix-valued coefficients, involving the generalized derivative that encompasses the Hukuhara, BG- and PS-derivatives. In the present work, we continue the study of non-homogeneous linear SVDEs with variable matrix-valued coefficients. The organization of the paper is as follows. In section 2, we present the basic notions required in our work and include other relevant results from [7, 10, 18]. Section 3 deals with the study of initial value problems associated with

---

[*1] Corresponding Author; e-mail: uepuganti2021@my.fit.edu

[2] e-mail: gtenali@fit.edu

non-homogeneous linear SVDEs with variable matrix-valued coefficients. Several illustrative examples are included.

## 2 Main Definitions and Notations

Let $K_c(\mathbb{R}^n)$ be a space of nonempty compact convex subsets of the Euclidean space $\mathbb{R}^n$ with the Hausdorff metric
$$h(A, B) = \min \{r \geq 0 : A \subset B + B_r(\mathbf{0}), B \subset A + B_r(\mathbf{0})\}$$
where $A, B \in K_c(\mathbb{R}^n)$ and $B_r(\mathbf{0}) = \{x \in \mathbb{R}^n : \|x\| \leq r\}$ is a closed ball of radius $r \geq 0$ centered at $\mathbf{0} \in \mathbb{R}^n$.

For any $A, B \in K_c(\mathbb{R}^n)$, and $\lambda \in \mathbb{R}$, the Minkowski addition and (scalar) multiplication are defined as below:
$$A + B = \{a + b : a \in A, b \in B\} \text{ and } \lambda A = \{\lambda a : a \in A, \lambda \in \mathbb{R}\}.$$

The following properties of these operations are well known and will be used throughout this paper. (see [3, 19]).

**Proposition 2.1.** For $A, B, C, D \in K_c(\mathbb{R}^n), \lambda, \mu \in \mathbb{R}_+$ and $\beta \in \mathbb{R}$:

(i) $A + \{\mathbf{0}\} = A$;
(ii) $A + B = B + A$;
(iii) $1A = A, 0A = \{\mathbf{0}\}$;
(iv) $(\lambda + \mu)A = \lambda A + \mu A$;
(v) $\beta(A + B) = \beta A + \beta B$;
(vi) $h(A + C, B + C) = h(A, B)$;
(vii) $h(\beta A, \beta B) = |\beta| h(A, B)$.

$(K_c(\mathbb{R}^n), h)$ is a complete metric space. For $A \in K_c(\mathbb{R}^n)$, in general, $A + (-1)A \neq \{\mathbf{0}\}$. To avoid this inconvenience, the notion of Hukuhara difference [2] is used.

**Definition 2.1.** ( [2]) For $A, B \in K_c(\mathbb{R}^n)$, if there exists an element $C \in K_c(\mathbb{R}^n)$ such that $A = B + C$, then the Hukuhara difference of $A$ and $B$ is said to exist and is denoted by $A \stackrel{h}{-} B$. i.e, $A \stackrel{h}{-} B = C \iff A = B + C$.

We list, below, some properties of the Hukuhara difference [3, 6, 19] that are used in the rest of the paper.

**Proposition 2.2.** ( [8])

(i) If the Hukuhara difference $A \stackrel{h}{-} B$ of two sets $A, B \in K_c(\mathbb{R}^n)$ exists, then it is unique;

(ii) $A \stackrel{h}{-} A = \{\mathbf{0}\}$ for all $A \in K_c(\mathbb{R}^n)$;

(iii) If $A \stackrel{h}{-} B$ exits and $\lambda \in \mathbb{R}_+$, then $\lambda A \stackrel{h}{-} \lambda B$ exits and $\lambda(A \stackrel{h}{-} B) = \lambda A \stackrel{h}{-} \lambda B$;

(iv) $(A + B) \stackrel{h}{-} B = A$ for all $A, B \in K_c(\mathbb{R}^n)$;

(v) If the Hukuhara difference $A \stackrel{h}{-} B$ of two sets $A, B \in K_c(\mathbb{R}^n)$ exists, then one can find $\alpha_0 \geq 1$ such that the Hukuhara difference $A \stackrel{h}{-} \alpha B$ exists for all $\alpha \in (0, \alpha_0]$;

(vi) If the Hukuhara difference $A \stackrel{h}{-} (B + C)$ exists, then the Hukuhara differences $A \stackrel{h}{-} B$ and $A \stackrel{h}{-} C$ also exist and, moreover, $A \stackrel{h}{-} (B + C) = A \stackrel{h}{-} B \stackrel{h}{-} C = A \stackrel{h}{-} C \stackrel{h}{-} B$;

(vii) If the Hukuhara differences $A \stackrel{h}{-} B$ and $A \stackrel{h}{-} C$ exist, then there is $\alpha > 0$ such that the Hukuhara difference $A \stackrel{h}{-} \alpha(B + C)$ exists.

**Definition 2.2.** ( [20])A set $M \subseteq \mathbb{R}^n$ is said to be *centrally symmetric* with respect to an element $q \in \mathbb{R}^n$ if for every $x \in M$, there exists $x' \in M$ such that $x + x' = 2q$.

We have the following additional properties involving Hukuhara differences of centrally symmetric sets.

**Lemma 2.1.** ( [10, 19])

(i) If the sets $A$ and $B$ are centrally symmetric and the Hukuhara difference $A \stackrel{h}{-} B$ exists, then the Hukuhara differences $A \stackrel{h}{-} (-1)A, B \stackrel{h}{-} (-1)B, (-1)A \stackrel{h}{-} B, A \stackrel{h}{-} (-1)B$, and $(-1)A \stackrel{h}{-} (-1)B$ exist;

(ii) If the sets $A$ and $B$ are not centrally symmetric and the Hukuhara difference $A \stackrel{h}{-} B$ exists, then the Hukuhara difference $(-1)A \stackrel{h}{-} (-1)B$ exists, whereas the Hukuhara differences $A \stackrel{h}{-} (-1)A, B \stackrel{h}{-} (-1)B, (-1)A \stackrel{h}{-} B$, and $A \stackrel{h}{-} (-1)B$ do not exist;

(iii) If $A, B \in K_c(\mathbb{R}^n)$ and the Hukuhara difference $A \stackrel{h}{-} B$ exists, then either both sets $A$ and $B$ are centrally symmetric or neither is centrally symmetric;

(iv) For $X, Y \in K_c(\mathbb{R}^n)$, if $X + Y = B_1(\mathbf{0})$, then $X = B_\mu(z_1)$ and $Y = B_\lambda(z_2)$ are such that $\mu + \lambda = 1$ and $z_1 + z_2 = \mathbf{0}$;

(v) $B_{r_1}(z_1) + B_{r_2}(z_2) = B_{r_1+r_2}(z_1 + z_2)$;

(vi) $\lambda B_r(\mathbf{0}) = B_{|\lambda|r}(\mathbf{0})$.

**Remark 2.1.** [7] If the difference $B_R(a) \stackrel{h}{-} X$ exists, then the set $X$ is also a ball $B_r(b)$ with radius $r$, where $r \leq R$.

Let $\mathcal{A} \in L(\mathbb{R}^n)$, the space of all linear operators on $\mathbb{R}^n$. The action of the linear operator $\mathcal{A}$ can be naturally extended to the space $K_c(\mathbb{R}^n)$ (the space of all compact convex subsets of $\mathbb{R}^n$) as follows:

$$\mathcal{A}X = \{\mathcal{A}x \mid x \in X\} \in K_c(\mathbb{R}^n). \tag{2.1}$$

**Proposition 2.3.** ( [19]) If $F, G \in K_c(\mathbb{R}^n)$, and $\mathcal{A}, \mathcal{B} \in L(\mathbb{R}^n)$, then the following relations are true:

(i) $\mathcal{A}(F + G) = \mathcal{A}F + \mathcal{A}G$;

(ii) $\mathcal{A}(\mathcal{B}F) = (\mathcal{A}\mathcal{B})F$;

(iii) $(\mathcal{A} + \mathcal{B})F \subseteq \mathcal{A}F + \mathcal{B}F$;

(iv) if $F \subseteq G$ then $\mathcal{A}F \subseteq \mathcal{A}F$.

We may consider for $X \in K_c(\mathbb{R}^n)$ and a $n \times n$ matrix $A$, the 'product' $AX$ as: $AX = \{Ax : x \in X\}$.

**Remark 2.2.** For $X = B_1(\mathbf{0})$, $AB_1(\mathbf{0}) = \{Ax : x \in B_1(\mathbf{0})\}$ is an $r$-dimensional ellipsoid and its semi-axes are equal to the corresponding singular values of the matrix $A$, where $r = \text{rank}(A)$ (see [21]). Further, It is clear that if the $n \times n$ matrix $A$ is orthogonal, then $AB_1(\mathbf{0}) = B_1(\mathbf{0})$.

**Definition 2.3.** ( [2]) Let $X(\cdot) : [0, T] \to K_c(\mathbb{R}^n)$ be the set-valued mapping. The Hukuhara derivative (H-derivative) $D_H X(t) \in K_c(\mathbb{R}^n)$ at $t \in (0, T)$ is defined if, for all sufficiently small $\Delta > 0$, the corresponding H-differences exist such that

$$\lim_{\Delta \to 0_+} \Delta^{-1}(X(t + \Delta) \stackrel{h}{-} X(t)) = \lim_{\Delta \to 0_+} \Delta^{-1}(X(t) \stackrel{h}{-} X(t - \Delta)) = D_H X(t).$$

**Remark 2.3.** If the mapping $X(\cdot)$ is H-differentiable on $[0, T]$, then the $\text{diam}(X(\cdot))$ is nondecreasing on $[0, T]$. The converse, however, is not true.

To address the shortcomings of the Hukuhara derivative, alternative notions of derivatives for set-valued mappings were introduced, such as the Huygens derivative, the $\pi$-derivative, and the T-derivative [22, 23]. These derivatives have been extensively studied, with their properties analyzed in several works [19, 24]. However, differential equations based on these derivatives tend to be highly complex, even at the formulation stage [19, 22]. Building on these advancements, the following two generalizations, namely, the BG-derivative (Bede and Gal) and the PS-derivatives (Plotnikov and Skripnik) were quite successful in aiming to overcome the intrinsic challenges posed by the use of Hukuhara derivative, and are used throughout our work.

**Definition 2.4.** ( [4, 14]) Let $X : [0, T] \to K_c(\mathbb{R}^n)$. The BG-derivative $D_{bg}X(t) \in K_c(\mathbb{R}^n)$ at $t \in (0, T)$ is defined if, for sufficiently small $\Delta > 0$, the corresponding H-differences exist and at least one of the following equalities holds:

(i) $\lim_{\Delta \to 0} \Delta^{-1}(X(t+\Delta) \stackrel{h}{-} X(t)) = \lim_{\Delta \to 0} \Delta^{-1}(X(t) \stackrel{h}{-} X(t-\Delta)) = D_{bg}X(t)$,

(ii) $\lim_{\Delta \to 0}(-\Delta)^{-1}(X(t) \stackrel{h}{-} X(t+\Delta)) = \lim_{\Delta \to 0}(-\Delta)^{-1}(X(t-\Delta) \stackrel{h}{-} X(t)) = D_{bg}X(t)$,

(iii) $\lim_{\Delta \to 0} \Delta^{-1}(X(t+\Delta) \stackrel{h}{-} X(t)) = \lim_{\Delta \to 0}(-\Delta)^{-1}(X(t-\Delta) \stackrel{h}{-} X(t)) = D_{bg}X(t)$,

(iv) $\lim_{\Delta \to 0}(-\Delta)^{-1}(X(t) \stackrel{h}{-} X(t+\Delta)) = \lim_{\Delta \to 0} \Delta^{-1}(X(t) \stackrel{h}{-} X(t-\Delta)) = D_{bg}X(t)$.

**Definition 2.5.** ( [6, 9]) Let $X : [0, T] \to K_c(\mathbb{R}^n)$. The PS-derivative $D_{ps}X(t) \in K_c(\mathbb{R}^n)$ at $t \in (0, T)$ is defined if, for sufficiently small $\Delta > 0$, the corresponding H-differences exist and at least one of the following equalities holds:

(i) $\lim_{\Delta \to 0} \Delta^{-1}(X(t+\Delta) \stackrel{h}{-} X(t)) = \lim_{\Delta \to 0} \Delta^{-1}(X(t) \stackrel{h}{-} X(t-\Delta)) = D_{ps}X(t)$,

(ii) $\lim_{\Delta \to 0} \Delta^{-1}(X(t) \stackrel{h}{-} X(t+\Delta)) = \lim_{\Delta \to 0} \Delta^{-1}(X(t-\Delta) \stackrel{h}{-} X(t)) = D_{ps}X(t)$,

(iii) $\lim_{\Delta \to 0} \Delta^{-1}(X(t+\Delta) \stackrel{h}{-} X(t)) = \lim_{\Delta \to 0} \Delta^{-1}(X(t-\Delta) \stackrel{h}{-} X(t)) = D_{ps}X(t)$,

(iv) $\lim_{\Delta \to 0} \Delta^{-1}(X(t) \stackrel{h}{-} X(t+\Delta)) = \lim_{\Delta \to 0} \Delta^{-1}(X(t) \stackrel{h}{-} X(t-\Delta)) = D_{ps}X(t)$.

**Remark 2.4.** ( [8]) If a set-valued mapping $X(\cdot)$ is H-differentiable on $[0, T]$, then it is BG-differentiable and PS-differentiable on $[0, T]$, and we have:

$$D_H X(t) = D_{ps}X(t) = D_{bg}X(t) \text{ for all } t \in [0, T].$$

There exist set-valued mappings that are BG-differentiable and PS-differentiable but not H-differentiable. And there exist set-valued mappings $X(\cdot)$ such that $D_{bg}X(t) \neq D_{ps}X(t)$ for all $t$. For more details, see [6].

**Definition 2.6.** ( [2]) Let $X : I = [t_0, T] \to K_c(\mathbb{R}^n)$. The Hukuhara Integral of $X$ over $I$ is defined as

$$\int_{t_0}^{t} X(s)\, ds = \left\{ \int_{t_0}^{t} x(s)\, ds : x \text{ is a measurable selector of } X \right\}.$$

We list below some properties of the Hukuhara Integral:

If $X, Y : [t_0, T] \to K_c(\mathbb{R}^n)$ are integrable, we have

$$\int_{t_0}^{t_2} X(s)\, ds = \int_{t_0}^{t_1} X(s)\, ds + \int_{t_1}^{t_2} X(s)\, ds, \quad t_0 \leq t_1 \leq t_2 \leq T,$$

$$\int_{t_0}^{t} \lambda X(s)\, ds = \lambda \int_{t_0}^{t} X(s)\, ds, \quad \lambda \in \mathbb{R}, \quad t_0 \leq t \leq T,$$

and,

$$h\left( \int_{t_0}^{t} X(s)\, ds, \int_{t_0}^{t} Y(s)\, ds \right) \leq \int_{t_0}^{t} h\left( X(s), Y(s) \right) ds.$$

**Theorem 2.1.** ( [2, 6, 14]) Suppose that the set-valued mapping $X : [t_0, T] \to K_c(\mathbb{R}^n)$ is such that one of $D_H X(\cdot), D_{ps} X(\cdot), or\ D_{bg} X(\cdot)$ exist on $[t_0, T]$.

(i) For all $t \in [t_0, T]$, if $\operatorname{diam}(X(t))$ is nondecreasing, then:

$$X(t) = X(t_0) + \int_{t_0}^{t} DX(s)ds,$$

where $DX(t)$ denotes one of the derivatives, $D_H X(t)$, $D_{ps} X(t)$ or $D_{bg} X(t)$.

(ii) For all $t \in [t_0, T]$, if the function $\operatorname{diam}(X(t))$ is decreasing, then:

$$X(t) = X(t_0) \stackrel{h}{-} \int_{t_0}^{t} D_{ps} X(s)ds;$$

$$X(t) = X(t_0) \stackrel{h}{-} (-1)\int_{t_0}^{t} D_{bg} X(s)ds.$$

From the results and discussions in [8, 9], we can see that the solutions of linear SVDEs with PS- and BG-derivatives coincide when the initial set $X_0$ and the set $F$ are centrally symmetric and, if in addition, $X_0 = (-1)X_0$, $F = (-1)F$. If these conditions are not met, the solutions may be different.

The linear homogeneous SVDEs with matrix coefficients are studied in [10]. In such cases, while it is possible to obtain general formulas for the solutions of such problems, it is necessary to make some additional assumptions on the matrix coefficients in order to obtain explicit formulas that can be used to construct set valued solutions and to have an idea of the underlying geometry of the solutions. For example, it is shown in [10] that when all the singular values of the coefficient matrix $A$ are equal, say to $\sigma$, two solutions of the initial value problem

$$DX(t) = AX(t), \quad X(0) = B_1(\mathbf{0})$$

are given by $X(t) = e^{\sigma t} B_1(\mathbf{0}), t \geq 0$ and $X(t) = e^{-\sigma t} B_1(\mathbf{0}), t \geq 0$. And, if the coefficient matrix has at least two different singular values, then the IVP has a solution only when the generalized derivative coincides with the Hukuhara derivative. Further, an explicit expression to the solution was given under the assumption that the coefficient matrix $A$ is either symmetric or positive definite. It was also established (see [10]) via a constructive formula that the solution of the IVP is an ellipsoid.

Following the work of Plotnikov, Skripnik and others in [10], we studied the initial value problems associated with non-homogeneous linear SVDEs, with (constant) matrix valued coefficients:

$$DX(t) = AX(t) + B_{r(t)}(\mathbf{0}), \quad X(0) = B_1(\mathbf{0}) \tag{2.2}$$

and obtained the form of the basic solutions and mixed solutions in [18]. A brief summary of the results obtained in [18] is given below:

(i) The initial value problem (2.2) with BG- and PS-derivatives admits infinitely many solutions. A solution is called a *basic solution* if the diameter of the (set-valued) solution is a monotone function, whereas the solutions whose diameters are not monotonic functions of time are referred to as *mixed solutions*. Specifically, the diameter $\operatorname{diam}(X_1(t))$ of the first basic solution $X_1(t)$ (corresponding to option (i) in the definitions 2.4 and 2.5) is a nondecreasing function whereas the second basic solution $X_2(t)$ (corresponding to option (ii) in the definitions 2.4 and 2.5) has a decreasing $\operatorname{diam}(X_2(t))$.

(ii) In the special case where the singular values of the matrix $A \in \mathbb{R}^{n \times n}$ are equal, i.e., $\sigma_1 = \cdots = \sigma_n = \sigma$, the differential equation (2.2) admits both basic solutions given by:

$$X_1(t) = \left(e^{\sigma t} + \int_0^t e^{\sigma(t-s)} r(s)\, ds\right) B_1(\mathbf{0}), \quad t \geq 0,$$

$$X_2(t) = \left(e^{-\sigma t} - \int_0^t e^{-\sigma(t-s)} r(s)\, ds\right) B_1(\mathbf{0}), \quad 0 \leq t \leq t_0,$$

where $t_0 = \frac{1}{\sigma} \ln\left(1 + \frac{\sigma}{\|r\|}\right)$. However, in the general case when the singular values of $A$ differ, the equation (2.2) possesses only the first basic solution $X_1(t)$ i.e, the solution which coincides with the Hukuhara (the diameter of $X_1(t)$ is nondecreasing).

The linear homogeneous SVDEs with variable matrix-valued coefficients were investigated in [7] using the generalized derivative (PS-derivative). The authors established conditions for the existence of solutions and derived explicit analytic expressions for the shapes of cross sections at every time under the assumption that the coefficient matrix possesses equal singular values and the system satisfies the Lappo-Danilevskii condition.

In this paper, continuing the study initiated in [18], we obtain the basic and mixed solutions for initial value problems associated non-homogeneous linear SVDEs, with variable matrix-valued coefficients. We present several insightful, illustrative examples. For convenience, we limit our study to mappings with values in $K_c(\mathbb{R}^2)$.

## 3 Linear Non-homogeneous Set-Valued Differential Equations

Consider the set-valued linear non-homogeneous differential equation of the form

$$DX(t) = A(t)X(t) + B_{r(t)}(\mathbf{0}), \quad X(t_0) = B_1(\mathbf{0}) \tag{3.1}$$

where $DX(\cdot)$ is one of the derivatives $D_H X(\cdot), D_{ps} X(\cdot)$, or $D_{bg} X(\cdot)$ of the set-valued mapping $X : [t_0, T] \to K_c(\mathbb{R}^2)$, $A(t) \in \mathbb{R}^{2\times 2}$ is a nondegenerate matrix for all $t \geq t_0 \geq 0$, and, $r(\cdot) : \mathbb{R} \to \mathbb{R}_+ = [t_0, \infty)$ is bounded continuous function.

Let $B_{r(t)}(\mathbf{0}) = \{x \in \mathbb{R}^2 : \|x\| \leq r(t)\}$. and, $R(t) = \int_{t_0}^t r(s)ds$.

The set-valued integral equation (SVIE) locally equivalent to the differential equation (3.1), depending on the type of derivative involved, is given below. For $t \in (t_0, T)$, if the term $DX(t)$ on the left hand side of (3.1) is

(i) $D_H X(t)$, then the SVIE is

$$X(t) = B_1(\mathbf{0}) + \int_{t_0}^t A(s)X(s)ds + B_{R(t)}(\mathbf{0}), \tag{3.2}$$

(ii) $D_{ps} X(t)$ due to (ii) in the definition (2.5), then the SVIE is

$$X(t) + \int_{t_0}^t A(s)X(s)ds + B_{R(t)}(\mathbf{0}) = B_1(\mathbf{0}), \tag{3.3}$$

(iii) $D_{bg} X(t)$ due to (ii) in the definition (2.4), then the SVIE is

$$X(t) + (-1)\int_{t_0}^t A(s)X(s)ds + B_{R(t)}(\mathbf{0}) = B_1(\mathbf{0}). \tag{3.4}$$

**Definition 3.1.** A set-valued mapping $X : [t_0, T] \to K_c(\mathbb{R}^2)$ is called a solution of the differential equation (3.1) with Hukuhara, PS, or BG-derivatives if it is continuous and satisfies the differential equation (3.1) (or a correspondent equivalent integral equation) on $[t_0, t_1]$, for some $t_1 < T$.

Let

$$A(t) = \begin{pmatrix} a(t) & b(t) \\ c(t) & d(t) \end{pmatrix}$$

where $a(t), b(t), c(t)$, and $d(t)$ are continuous functions such that $\det(A(t)) \neq 0$ for all $t \geq t_0 \geq 0$.

For $t \geq t_0$, the singular values of the matrix $A(t)$ are given by:

$$\sigma_{1,2}(t) = \sqrt{\frac{(\|A(t)\|)^2 \pm \sqrt{\delta(t)}}{2}},$$

where $\|A(t)\|$ is Frobenius norm, and $\delta(t) = (\|A(t)\|)^2 - 4(\det A(t))^2 \geq 0$.

**Remark 3.1.** The singular values of $A(t)$ are equal if and only if $\delta(t) \equiv 0$ or equivalently if and only if $A(t)$ is either

$$A(t) = \begin{pmatrix} a(t) & b(t) \\ -b(t) & a(t) \end{pmatrix} \quad \text{or} \quad A(t) = \begin{pmatrix} a(t) & b(t) \\ b(t) & -a(t) \end{pmatrix}.$$

In this case, the singular values are given by

$$\sigma_1(t) = \sigma_2(t) = \sigma(t) = \frac{\|A(t)\|}{\sqrt{2}}.$$

In the following theorem, we obtain an explicit form for the (two) basic solutions of the initial value problem (3.1), when the coefficient matrix $A(t)$ has equal singular values.

**Theorem 3.1.** If the matrix $A(t)$ has two equal singular values, say, $\sigma_1(t) = \sigma_2(t) = \sigma(t)$, then the IVP (3.1) with BG- or PS-derivatives admits two basic solutions:

(i) The first basic solution $X_1(\cdot)$ (corresponding to Hukuhara derivative or option (i) in the definitions 2.4 and 2.5) is given by

$$X_1(t) = \left(G(t_0, t) + \int_{t_0}^{t} G(s, t) r(s) \, ds\right) B_1(\mathbf{0}), \quad G(s, t) = e^{\int_s^t \sigma(u) \, du}, \quad t \geq s \geq t_0 \geq 0.$$

(ii) The second basic solution $X_2(\cdot)$ (corresponding to option (ii) in the definitions 2.4 and 2.5) is given by

$$X_2(t) = \left(H(t_0, t) - \int_{t_0}^{t} H(s, t) r(s) \, ds\right) B_1(\mathbf{0}), \quad H(s, t) = e^{-\int_s^t \sigma(u) \, du}, \quad 0 \leq t_0 \leq s \leq t < \tau,$$

where $\tau > 0$ is such that $\int_{t_0}^{\tau} e^{\int_{t_0}^{s} \sigma(u) \, du} r(s) ds = 1$.

**Remark 3.2.** Clearly, the interval of existence of the second basic solution is dependent on the nature of $\|A(t)\|$. For example, the solution exists on a larger interval if $\|A(t)\|$ is a monotonically decreasing function and on a smaller interval if $\|A(t)\|$ is a monotonically increasing function. Further, if $r(t) = \|A(t)\|$, then the interval of the existence of the second basic solution is $[t_0, \tau)$, where $\tau$ is determined by $\int_{t_0}^{\tau} \|A(u)\| du = \sqrt{2} \ln\left(1 + \frac{1}{\sqrt{2}}\right)$.

**Proof:** First, we prove that $X_1(\cdot)$ is a first basic solution of the IVP (3.1). To this end, we demonstrate that the following identity holds:

$$D\left(\left(G(t_0, t) + \int_{t_0}^{t} G(s, t) r(s) \, ds\right) B_1(\mathbf{0})\right) \equiv A(t) \left(G(t_0, t) + \int_{t_0}^{t} G(s, t) r(s) \, ds\right) B_1(\mathbf{0}) + B_{r(t)}(\mathbf{0}).$$

Since $\sigma(t) > 0$ for all $t \geq t_0$, we conclude that $\rho(t) = G(t_0, t) + \int_{t_0}^{t} G(s, t) r(s) \, ds$ is a nondecreasing function and, hence, the same is true for $\text{diam}(X_1(\cdot))$. It is clear that

$$G(t_0, t + \Delta) - G(t_0, t) = G(t_0, t) \left(G(t, t + \Delta) - 1\right),$$

and,

$$\int_{t_0}^{t+\Delta} G(s, t+\Delta) r(s) \, ds - \int_{t_0}^{t} G(s, t) r(s) \, ds = \int_{t_0}^{t} G(s, t) r(s) \, ds \left(G(t, t+\Delta) - 1\right) + \int_{t}^{t+\Delta} G(s, t+\Delta) r(s) \, ds.$$

Thus, according to the Definition (2.3), and using the L'Hôpital's rule, we have,

$$\lim_{\Delta \to 0_+} \Delta^{-1}\left(X_1(t+\Delta) \stackrel{h}{-} X_1(t)\right)$$

$$= \lim_{\Delta \to 0_+} \Delta^{-1}\left(G(t_0, t+\Delta) - G(t_0, t) + \int_{t_0}^{t+\Delta} G(s, t+\Delta)\, r(s)\, \mathrm{d}s - \int_{t_0}^{t} G(s, t)\, r(s)\, \mathrm{d}s\right) B_1(\mathbf{0})$$

$$= \lim_{\Delta \to 0_+} \Delta^{-1}\left(\left(G(t_0, t) + \int_{t_0}^{t} G(s, t)\, r(s)\, \mathrm{d}s\right)\left(G(t, t+\Delta) - 1\right) + \int_{t}^{t+\Delta} G(s, t+\Delta)\, r(s)\, \mathrm{d}s\right) B_1(\mathbf{0})$$

$$= \lim_{\Delta \to 0_+} \left(\left(G(t_0, t) + \int_{t_0}^{t} G(s, t)\, r(s)\, \mathrm{d}s\right) \sigma(t+\Delta) G(t, t+\Delta) + G(t+\Delta, t+\Delta) r(t+\Delta)\right) B_1(\mathbf{0})$$

$$= \left(\sigma(t)\left(G(t_0, t) + \int_{t_0}^{t} G(s, t)\, r(s)\, \mathrm{d}s\right) + r(t)\right) B_1(\mathbf{0}).$$

Similarly, it is clear that

$$G(t_0, t) - G(t_0, t-\Delta) = G(t_0, t)\left(1 - G(t, t-\Delta)\right),$$

and,

$$\int_{t_0}^{t} G(s, t)\, r(s)\, \mathrm{d}s - \int_{t_0}^{t-\Delta} G(s, t-\Delta)\, r(s)\, \mathrm{d}s = \int_{t_0}^{t} G(s, t)\, r(s)\, \mathrm{d}s\left(1 - G(t, t-\Delta)\right) + \int_{t-\Delta}^{t} G(s, t-\Delta)\, r(s)\, \mathrm{d}s.$$

Thus,

$$\lim_{\Delta \to 0_+} \Delta^{-1}\left(X_1(t) \stackrel{h}{-} X_1(t-\Delta)\right)$$

$$= \lim_{\Delta \to 0_+} \Delta^{-1}\left(G(t_0, t) - G(t_0, t-\Delta) + \int_{t_0}^{t} G(s, t)\, r(s)\, \mathrm{d}s - \int_{t_0}^{t-\Delta} G(s, t-\Delta)\, r(s)\, \mathrm{d}s\right) B_1(\mathbf{0})$$

$$= \lim_{\Delta \to 0_+} \Delta^{-1}\left(\left(G(t_0, t) + \int_{t_0}^{t} G(s, t)\, r(s)\, \mathrm{d}s\right)\left(1 - G(t, t-\Delta)\right) + \int_{t-\Delta}^{t} G(s, t-\Delta)\, r(s)\, \mathrm{d}s\right) B_1(\mathbf{0})$$

$$= \lim_{\Delta \to 0_+} \left(\left(G(t_0, t) + \int_{t_0}^{t} G(s, t)\, r(s)\, \mathrm{d}s\right)\left(\sigma(t-\Delta) G(t, t-\Delta)\right) + G(t-\Delta, t-\Delta) r(t-\Delta)\right) B_1(\mathbf{0})$$

$$= \left(\sigma(t)\left(G(t_0, t) + \int_{t_0}^{t} G(s, t)\, r(s)\, \mathrm{d}s\right) + r(t)\right) B_1(\mathbf{0}).$$

Therefore,

$$DX_1(t) = \left(\sigma(t)\left(G(t_0, t) + \int_{t_0}^{t} G(s, t)\, r(s)\, \mathrm{d}s\right) + r(t)\right) B_1(\mathbf{0})$$

Since the matrix $A(t)$ has equal singular values, i.e., $\sigma_1(t) = \sigma_2(t) = \sigma(t)$, its singular value decomposition at any time $t \geq t_0$ takes the form

$$A(t) = U(t)\Sigma(t)V^T(t),$$

where $U(t)$ and $V(t)$ are orthogonal matrices, and $\Sigma(t) = \sigma(t)I$, where $I$ is the identity matrix. It is also known that $V^T(t)B_r(\mathbf{0}) = B_r(\mathbf{0})$ and $U(t)B_r(\mathbf{0}) = B_r(\mathbf{0})$ for all $t \geq t_0$ and $r \geq 0$.

And we have,

$$A(t)B_1(\mathbf{0}) = U(t)\Sigma(t)V^T(t)B_1(\mathbf{0}) = \sigma(t)U(t)V^T(t)B_1(\mathbf{0}) = \sigma(t)U(t)B_1(\mathbf{0}) = \sigma(t)B_1(\mathbf{0})$$

Now consider the RHS of the equation (3.1),

$$A(t)\left(G(t_0,t) + \int_{t_0}^{t} G(s,t)\,r(s)\,\mathrm{d}s\right) B_1(\mathbf{0}) + B_{r(t)}(\mathbf{0}) = A(t)\rho(t) B_1(\mathbf{0}) + B_{r(t)}(\mathbf{0})$$
$$= \sigma(t) U(t) V^T(t) B_{\rho(t)}(\mathbf{0}) + B_{r(t)}(\mathbf{0})$$
$$= \left(\sigma(t)\left(G(t_0,t) + \int_{t_0}^{t} G(s,t)\,r(s)\,\mathrm{d}s\right) + r(t)\right) B_1(\mathbf{0}).$$

Therefore, it follows that $X_1(t)$ is the first basic solution of the IVP (3.1).

Since both the initial set and forcing terms are centrally symmetric and are such that $X(0) = (-1)X(0)$, and $F = (-1)F$, the second basic solution to the IVP (3.1) with PS-derivative coincides with the second basic solution to the IVP (3.1) with BG-derivative. Thus, it is sufficient to prove that the second basic solution is a solution of the IVP (3.1) with the PS-derivative; the corresponding result for the BG-derivative case then follows immediately.

Now, to establish that $X_2(\cdot)$ is a second basic solution of the IVP (3.1), we verify that the below identity holds with the generalized derivative corresponding to definition (2.5)(ii) of the PS-derivative.

$$D\left(\left(H(t_0,t) - \int_{t_0}^{t} H(s,t)\,r(s)\,\mathrm{d}s\right) B_1(\mathbf{0})\right) \equiv A(t)\left(H(t_0,t) - \int_{t_0}^{t} H(s,t)\,r(s)\,\mathrm{d}s\right) B_1(\mathbf{0}) + B_{r(t)}(\mathbf{0}).$$

Since $\sigma(t) > 0$ for all $t \geq t_0$, we conclude that $\omega(t) = H(t_0,t) - \int_{t_0}^{t} H(s,t)\,r(s)\,\mathrm{d}s$ is a decreasing function on $[t_0, \tau)$, where $\tau > 0$ is such that $\int_{t_0}^{\tau} e^{\int_{t_0}^{s} \sigma(u)\,du} r(s)\,ds = 1$.

Similarly, it is clear that

$$H(t_0,t) - H(t_0, t+\Delta) = H(t_0,t)\left(1 - H(t, t+\Delta)\right),$$

and,

$$\int_{t_0}^{t+\Delta} H(s, t+\Delta)\,r(s)\,\mathrm{d}s - \int_{t_0}^{t} H(s,t)\,r(s)\,\mathrm{d}s = \int_{t}^{t+\Delta} H(s, t+\Delta)\,r(s)\,\mathrm{d}s - \int_{t_0}^{t} H(s,t)\,r(s)\,\mathrm{d}s\left(1 - H(t, t+\Delta)\right).$$

Since the set $B_1(\mathbf{0})$ is centrally symmetric, and $(-1)B_1(\mathbf{0}) = B_1(\mathbf{0})$, using the L'Hôpital's rule we get,

$$\lim_{\Delta \to 0_+} \Delta^{-1}\left(X_2(t) \stackrel{h}{-} X_2(t+\Delta)\right)$$
$$= \lim_{\Delta \to 0_+} \Delta^{-1}\left(H(t_0,t) - H(t_0, t+\Delta) - \int_{t_0}^{t} H(s,t)\,r(s)\,\mathrm{d}s + \int_{t_0}^{t+\Delta} H(s, t+\Delta)\,r(s)\,\mathrm{d}s\right) B_1(\mathbf{0})$$
$$= \lim_{\Delta \to 0_+} \Delta^{-1}\left(\left(H(t_0,t) - \int_{t_0}^{t} H(s,t)\,r(s)\,\mathrm{d}s\right)\left(1 - H(t, t+\Delta)\right) + \int_{t}^{t+\Delta} H(s, t+\Delta)\,r(s)\,\mathrm{d}s\right) B_1(\mathbf{0})$$
$$= \lim_{\Delta \to 0_+}\left(\left(H(t_0,t) - \int_{t_0}^{t} H(s,t)\,r(s)\,\mathrm{d}s\right)\sigma(t+\Delta) H(t, t+\Delta) + H(t+\Delta, t+\Delta) r(t+\Delta)\right) B_1(\mathbf{0})$$
$$= \left(\sigma(t)\left(H(t_0,t) - \int_{t_0}^{t} H(s,t)\,r(s)\,\mathrm{d}s\right) + r(t)\right) B_1(\mathbf{0}).$$

It is clear that

$$H(t_0, t-\Delta) - H(t_0,t) = H(t_0,t)\left(H(t, t-\Delta) - 1\right),$$

and,

$$\int_{t_0}^{t} H(s,t)\,r(s)\,\mathrm{d}s - \int_{t_0}^{t-\Delta} H(s, t-\Delta)\,r(s)\,\mathrm{d}s = \int_{t-\Delta}^{t} H(s, t-\Delta)\,r(s)\,\mathrm{d}s - \int_{t_0}^{t} H(s,t)\,r(s)\,\mathrm{d}s\left(H(t, t-\Delta) - 1\right).$$

Thus, we have,

$$\lim_{\Delta \to 0_+} \Delta^{-1}\left(X_2(t-\Delta) \overset{h}{-} X_2(t)\right)$$

$$= \left(H(t_0, t-\Delta) - H(t_0, t) - \int_{t_0}^{t-\Delta} H(s, t-\Delta)\,r(s)\,ds + \int_{t_0}^{t} H(s,t)\,r(s)\,ds\right) B_1(\mathbf{0})$$

$$= \left(\left(H(t_0, t) - \int_{t_0}^{t} H(s,t)\,r(s)\,ds\right)\left(H(t, t-\Delta) - 1\right) + \int_{t-\Delta}^{t} H(s, t-\Delta)\,r(s)\,ds\right) B_1(\mathbf{0})$$

$$= \lim_{\Delta \to 0_+} \left(\left(H(t_0, t) - \int_{t_0}^{t} H(s,t)\,r(s)\,ds\right) H(t, t-\Delta)\sigma(t-\Delta) + H(t-\Delta, t-\Delta) r(t-\Delta)\right) B_1(\mathbf{0})$$

$$= \left(\sigma(t)\left(H(t_0, t) - \int_{t_0}^{t} H(s,t)\,r(s)\,ds\right) + r(t)\right) B_1(\mathbf{0}).$$

Therefore,

$$D_{ps}X_2(t) = \left(\sigma(t)\left(H(t_0, t) - \int_{t_0}^{t} H(s,t)\,r(s)\,ds\right) + r(t)\right) B_1(\mathbf{0})$$

Now consider the RHS of the equation (3.1),

$$A(t)\left(H(t_0, t) - \int_{t_0}^{t} H(s,t)\,r(s)\,ds\right) B_1(\mathbf{0}) + B_{r(t)}(\mathbf{0}) = A(t)\omega(t)B_1(\mathbf{0}) + B_{r(t)}(\mathbf{0})$$

$$= \sigma(t) U(t) V^T(t) B_{\omega(t)}(\mathbf{0}) + B_{r(t)}(\mathbf{0})$$

$$= \left(\sigma(t)\left(H(t_0, t) - \int_{t_0}^{t} H(s,t)\,r(s)\,ds\right) + r(t)\right) B_1(\mathbf{0})$$

Thus it follows that $X_2(t)$ is the second basic solution of the IVP (3.1).

**Example 3.1.** Consider the initial value problem

$$DX(t) = A(t)X(t) + B_{r(t)}(\mathbf{0}), \quad X(1) = B_1(\mathbf{0})$$

where $DX(t)$ is one of the derivatives $D_H X(t), D_{ps}X(t)$, or $D_{bg}X(t)$ of the set-valued map $X(t): [1, T] \to K_c(\mathbb{R}^2)$, and

$$A(t) = \begin{pmatrix} \frac{1}{t}\cos(t) & \frac{1}{t}\sin(t) \\ -\frac{1}{t}\sin(t) & \frac{1}{t}\cos(t) \end{pmatrix}$$

Then the singular values $\sigma_1(t)$ and $\sigma_2(t)$ of the matrix $A(t)$ are $\sigma_1(t) \equiv \sigma_2(t) = \sigma(t) = \frac{1}{t}$.

(i) If $r(t) = t$, then the first and second basic solutions to the IVP (3.1) are

$$X_1(t) = t^2 B_1(\mathbf{0}), \quad \forall t \in [1, +\infty), \quad \text{(Fig. (1))},$$

and,

$$X_2(t) = \left(\frac{4 - t^3}{3t}\right) B_1(\mathbf{0}), \quad \forall t \in [1, 1.5874), \quad \text{(Fig. (2))}.$$

(ii) If $r(t) = \frac{1}{t}$, then the first and second basic solutions to the IVP (3.1) are

$$X_1(t) = (2t - 1) B_1(\mathbf{0}), \quad \forall t \in [1, +\infty), \quad \text{(Fig. (3))},$$

and,

$$X_2(t) = \left(\frac{2}{t} - 1\right) B_1(\mathbf{0}), \quad \forall t \in [1, 2), \quad \text{(Fig. (4))}.$$

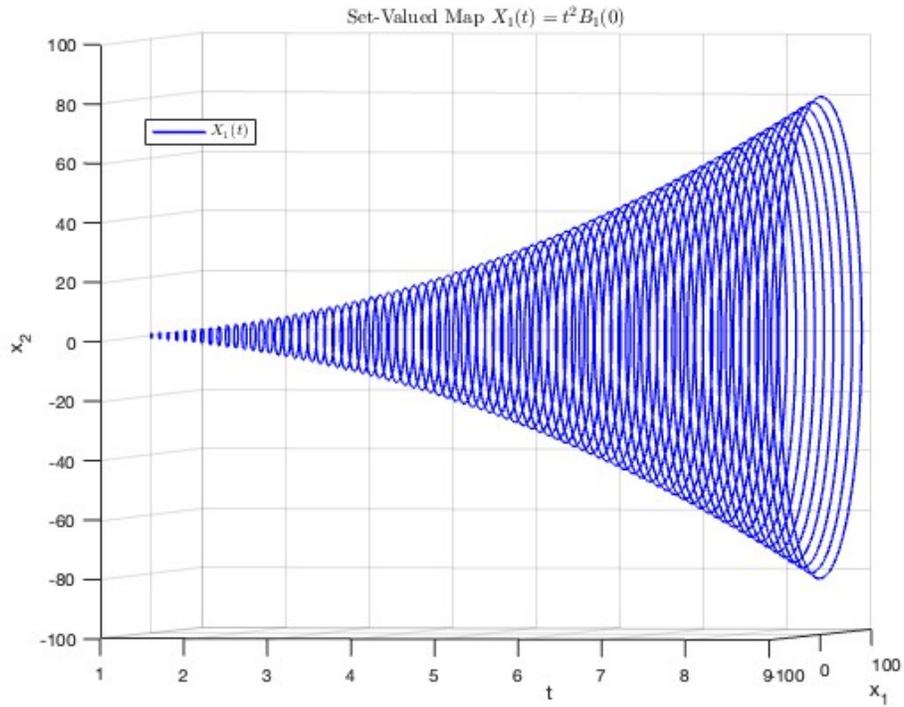

Figure 1: Solution $X_1(t), t \in [1, +\infty)$

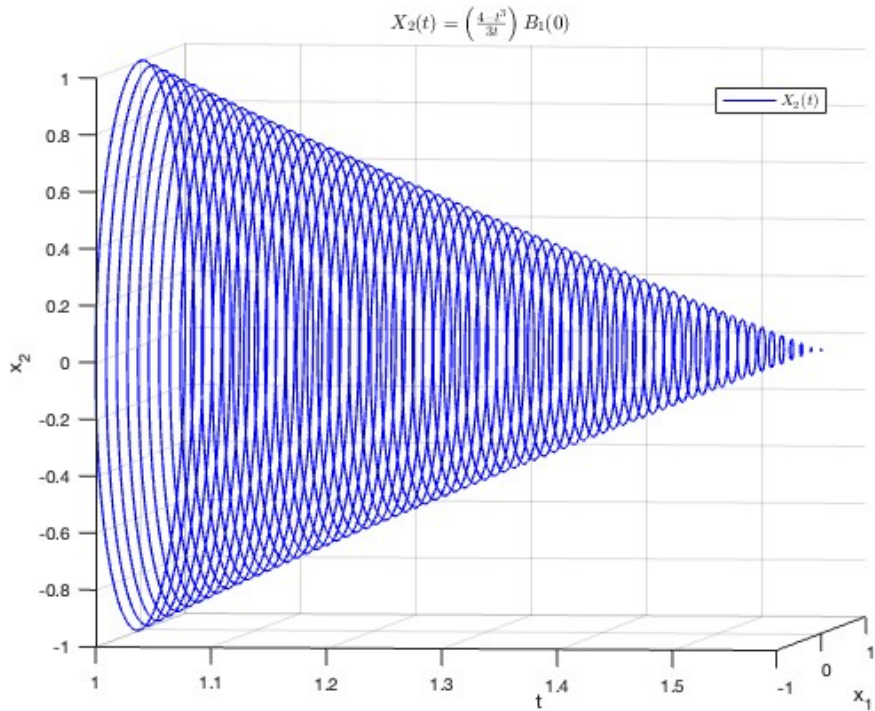

Figure 2: Solution $X_2(t), t \in [1, 1.5874)$

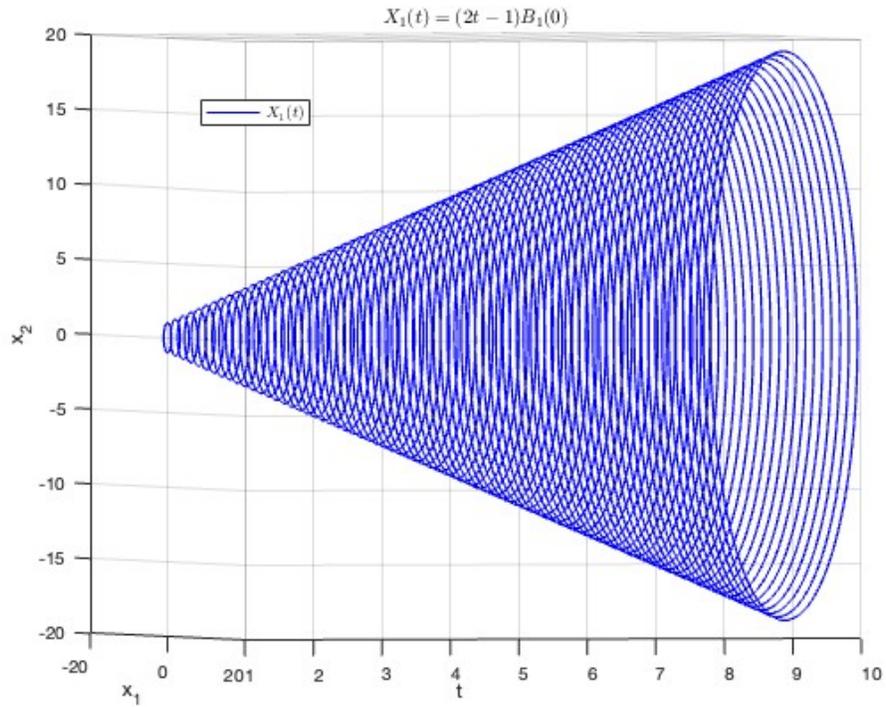

Figure 3: Solution $X_1(t), t \in [1, +\infty)$

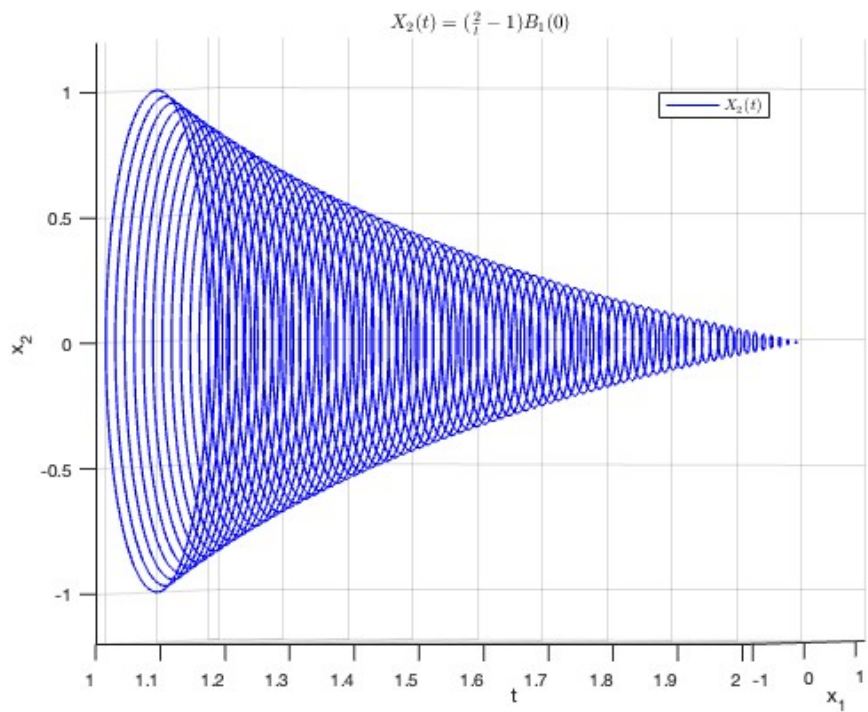

Figure 4: Solution $X_2(t), t \in [1, 2)$

Now, we consider the case where the matrix $A(t)$ is such that $\delta(t) \neq 0$, i.e, the singular values of $A(t)$ are not equal.

**Theorem 3.2.** If the matrix $A(t)$ is such that

$$A(t) = \begin{pmatrix} a(t) & b(t) \\ c(t) & d(t) \end{pmatrix}$$

and there exists $t' > t_0 \geq 0$ for which $\delta(t') \neq 0$, then the equation (3.1) with BG- or PS-derivatives has only the first basic solution $X_1(t)$, i.e the solution which coincides with the solution to the equation (3.1) with Hukuhara derivative.

**Proof:** It is known that the first basic solution of equation (3.1) coincides with the solution of the corresponding differential equation involving the Hukuhara derivative. We prove by contradiction that the equation (3.1) does not possess a second basic solution $X_2(\cdot)$. Assume, on the contrary that the second basic solution of the equation (3.1), $X_2(\cdot)$ exists.

Since there exists $t' > t_0 \geq 0$ such that $\delta(t') \neq 0$ and the function $\delta(t)$ is continuous on $\mathbb{R}_+$, there exists an interval $(t_1, t_2)$ such that
(i) $t_1 > t_0$ and $\delta(t_1) = 0$ or $t_1 = t_0$;
(ii) $t' \in (t_1, t_2)$;
(iii) $\delta(t) > 0$ for all $t \in (t_1, t_2)$.

For some $t_1 > t_0$, suppose that $\delta(t) = 0$ for all $t \in [t_0, t_1]$. (Note that if this is not the case, then there exists $t'' \in (t_0, t_1)$ such that $\delta(t'') > 0$, in which case, we replace $t'$ by $t''$ and proceed analogously.) Since $\delta(t) = 0$ for all $t \in [t_0, t_1]$, we get $\sigma_1(t) = \sigma_2(t) = \sigma(t)$ and, by virtue of Theorem (3.1),

$$X_2(t_1) = \left( H(t_0, t_1) - \int_{t_0}^{t_1} H(s, t_1) r(s) \, ds \right) B_1(\mathbf{0}) = B_{\omega(t_1)}(\mathbf{0}), \quad t_1 \geq s \geq t_0 \geq 0.$$

where $\omega(t_1) = \left( H(t_0, t_1) - \int_{t_0}^{t_1} H(s, t_1) r(s) \, ds \right) < 1$. If $t_1 = t_0$, then $X_2(t_1) = B_1(\mathbf{0})$.

By virtue of Definition 2.1 and Theorem 2.1, for all $t \in [t_1, t_2]$,

$$X_2(t) + \int_{t_1}^{t} A(s) X_2(s) ds + B_{R(t)}(\mathbf{0}) = B_{\omega(t_1)}(\mathbf{0}).$$

For a fixed but an arbitrary $T \in (t_1, t_2]$,

$$X_2(T) + \int_{t_1}^{T} A(s) X_2(s) ds + B_{R(T)}(\mathbf{0}) = B_{\omega(t_1)}(\mathbf{0}).$$

Thus,

$$X_2(T) + \int_{t_1}^{T} A(s) X_2(s) ds = B_{\omega(t_1)}(\mathbf{0}) \stackrel{h}{-} B_{R(T)}(\mathbf{0}). \tag{3.5}$$

This yields

$$B_{\omega(t_1) - R(T)}(\mathbf{0}) \stackrel{h}{-} X_2(T) = \int_{t_1}^{T} A(s) X_2(s) ds.$$

Since the Hukuhara difference $B_{\omega(t_1) - R(T)}(\mathbf{0}) \stackrel{h}{-} X_2(T)$ exists, by Lemma (2.1)(iv) and Remark (2.1), it follows that $X_2(T)$ is also a disk, i.e., $X_2(T) = B_{\rho(T)}(\mathbf{0})$, where $0 \leq \rho(T) \leq \omega(t_1) - R(T)$. Since $T$ is arbitrary, we have $X_2(t) = B_{\rho(t)}(\mathbf{0})$ for all $t \in [t_1, t_2]$.

Since the matrix $A(t)$ is such that $\delta(t) \neq 0$ for all $t \in (t_1, t_2]$, we conclude that, the matrix $A(t)$ has singular values $\sigma_1 \neq \sigma_2$, and $A(t)B_{\rho(t)}(\mathbf{0})$ is an ellipse with semi axes $e^{\sigma_1 t}$ and $e^{\sigma_2 t}$, for $t \in (t_1, t_2]$. Therefore,

$$\int_{t_1}^{T} A(s)X_2(s)ds = \int_{t_1}^{T} A(s)B_{\rho(s)}(\mathbf{0})ds$$

cannot be a disk. Hence,

$$B_{\rho(T)}(\mathbf{0}) + \int_{t_1}^{T} A(s)X_2(s)ds$$

is not a disk, and equality (3.5) is not true.

Therefore,

(i) if $\delta(t) > 0$ for all $t \in (t_0, t_2)$, then the second basic solution of the differential equation (3.1) does not exist for any $t \geq t_0$;

(ii) if $\delta(t) = 0$ for all $t \in [t_0, t_1]$ and $\delta(t) > 0$ for all $t \in (t_1, t_2)$, then the second basic solution of the differential equation (3.1) exists on $t \in [t_0, t_1]$ but does not exist for all $t > t_1$.

Hence, according to Definition (3.1), the second basic solution of the differential equation (3.1) does not exist. Theorem (3.2) is proved.

Even though the coefficient matrix $A(t)$ does not have equal singular values, it is possible to obtain an explicit form for the first basic solution of (3.1), when $A(t)$ satifies Lappo-Danilevskii condition and is symmetric.

**Theorem 3.3.** ( [25], p.24) A non-diagonal matrix $A(t)$ satisfies the Lappo-Danilevskii condition

$$A(t)\int_{t_0}^{t} A(s)ds = \left(\int_{t_0}^{t} A(s)\,ds\right) A(t)$$

for all $t \geq t_0 \geq 0$ if and only if it has the form:

$$A(t) = \begin{pmatrix} p(t) + \gamma q(t) & q(t) \\ \mu q(t) & p(t) \end{pmatrix}$$

where $p(t), q(t) : \mathbb{R}_+ \to \mathbb{R}$ are continuous functions and $\mu, \gamma \in \mathbb{R}$ are constants.

In what follows, we assume that the coefficient matrix $A(t) = \begin{pmatrix} p(t) + \gamma q(t) & q(t) \\ q(t) & p(t) \end{pmatrix}$, where $p(t), q(t)$ are continuous functions and $\gamma \in \mathbb{R}$ is a constant such that $q(t) \neq 0$ and $p(t) + \gamma q(t) \neq -p(t)$ for all $t \geq t_0 \geq 0$.

We now proceed to derive the explicit form of the first basic solution of (3.1).

Let

$$G(s,t) = \begin{pmatrix} \int_s^t (p(w) + \gamma q(w))\,dw & \int_s^t q(w)\,dw \\ \int_s^t q(w)\,dw & \int_s^t p(w)\,dw \end{pmatrix}.$$

Since the matrix $G(s,t)$ is symmetric for all $t \geq t_0$, at any time $t \geq t_0$, it has the following representation:

$$G(s,t) = U(s,t)\Lambda(s,t)U^T(s,t),$$

where $U(s,t)$ is an orthogonal matrix at any time $t \geq t_0$ is such that

$$U(s,t) = \begin{pmatrix} \frac{\eta(s,t)}{\sqrt{\eta^2(s,t)+4}} & \frac{2}{\sqrt{\eta^2(s,t)+4}} \\ \frac{2}{\sqrt{\eta^2(s,t)+4}} & \frac{-\eta(s,t)}{\sqrt{\eta^2(s,t)+4}} \end{pmatrix}, \tag{3.6}$$

and, $\Lambda(s,t) = \begin{pmatrix} \lambda_1(s,t) & 0 \\ 0 & \lambda_2(s,t) \end{pmatrix}$, $\lambda_i,(\cdot): \mathbb{R} \times \mathbb{R} \to \mathbb{R}, i=1,2, |\lambda_1(s,t)| > |\lambda_2(s,t)|$ is such that

$$\lambda_{1,2}(s,t) = \frac{1}{2}\left(\int_s^t (2p(w) + \gamma q(w))dw \pm \sqrt{(4+\gamma^2)}\left|\int_s^t q(w)dw\right|\right).$$

Here, $\eta(s,t) = \gamma + \sqrt{\gamma^2+4}\,\text{sgn}\left(\int_s^t q(w)\,dw\right)$, and $\text{sgn}(\varphi) = \begin{cases} 1, & \varphi \geq 0 \\ -1, & \varphi < 0 \end{cases}$.

Since $q(t) \neq 0$ for all $t \geq t_0$, we conclude that $q(t) > 0$ for all $t \geq t_0$ or $q(t) < 0$ for all $t \geq t_0$, and both $\int_{t_0}^t q(w)dw$ and $\int_s^t q(w)dw$ have the same sign. Therefore, $U(s,t) \equiv U(t_0,t) \equiv U$ for all $s,t \geq t_0 \geq 0$.

Further, using the notation: $\Sigma_1(t) = |\Lambda(t_0,t)|$, $\Sigma_2(t) = \int_{t_0}^t |\Lambda(s,t)|r(s)\,ds$, $|\Lambda(s,t)| = \Lambda(s,t)P(s,t)$ where the matrix

$$P(s,t) = \begin{pmatrix} \text{sgn}(\lambda_1(s,t)) & 0 \\ 0 & \text{sgn}(\lambda_2(s,t)) \end{pmatrix},$$

the first basic solution of the differential equation (3.1) can be realized in a more useful form.

**Theorem 3.4.** If the coefficient matrix

$$A(t) = \begin{pmatrix} p(t) + \gamma q(t) & q(t) \\ q(t) & p(t) \end{pmatrix}$$

where $p(t), q(t)$ are continuous functions such that $q(t) \neq 0$ and $p(t) + \gamma q(t) \neq -p(t)$ for all $t \geq t_0$, and $\gamma \in \mathbb{R}$, then the IVP (3.1) has only the first basic solution $X_1(t)$, given by

$$X_1(t) = U\left(e^{\Sigma_1(t)}B_1(\mathbf{0}) + e^{\Sigma_2(t)}B_1(\mathbf{0})\right).$$

**Proof:** Since $p(t), q(t)$, and $\gamma$ are such that $q(t) \neq 0$ and $p(t) + \gamma q(t) \neq -p(t)$ for all $t \geq t_0$, we get $\delta(t) \neq 0$ and, by Theorem (3.2), the IVP (3.1) has only the first basic solution $X_1(\cdot)$.

We represent the solution $X_1(\cdot)$ of the differential equation (3.1) in the form

$$X_1(t) = B_1(\mathbf{0}) + G(t_0,t)B_1(\mathbf{0}) + \frac{(G(t_0,t))^2}{2!}B_1(\mathbf{0}) + \ldots + \frac{(G(t_0,t))^k}{k!}B_1(\mathbf{0}) + \ldots + \int_{t_0}^t B_{r(s)}(\mathbf{0})\,ds$$

$$+ \int_0^t G(s,t)B_{r(s)}(\mathbf{0})\,ds + \int_{t_0}^t \frac{(G(s,t))^2}{2!}B_{r(s)}(\mathbf{0})\,ds + \ldots + \int_{t_0}^t \frac{(G(s,t))^k}{k!}B_{r(s)}(\mathbf{0})\,ds + \ldots$$

Since

$$\frac{(G(t_0,t))^k}{k!}B_1(\mathbf{0}) + \int_{t_0}^t \frac{(G(s,t))^k}{k!}B_{r(s)}(\mathbf{0})\,ds = \frac{1}{k!}U\Lambda^k(t_0,t)U^T B_1(\mathbf{0}) + \frac{1}{k!}\int_{t_0}^t U\Lambda^k(s,t)U^T B_{r(s)}(\mathbf{0})\,ds$$

$$= \frac{1}{k!}U\Lambda^k(t_0,t)B_1(\mathbf{0}) + \frac{1}{k!}\int_{t_0}^t U\Lambda^k(s,t)B_{r(s)}(\mathbf{0})\,ds$$

$$= \frac{1}{k!}U|\Lambda(t_0,t)|^k P^k(0,t)B_1(\mathbf{0}) + \frac{1}{k!}\int_{t_0}^t U|\Lambda^k(s,t)|P^k(s,t)B_{r(s)}(\mathbf{0})\,ds$$

$$= \frac{1}{k!}U|\Lambda(t_0,t)|^k B_1(\mathbf{0}) + \frac{1}{k!}U\int_{t_0}^t |\Lambda^k(s,t)|B_{r(s)}(\mathbf{0})\,ds$$

$$= \frac{1}{k!}U\Sigma_1^k(t)B_1(\mathbf{0}) + \frac{1}{k!}U\Sigma_2^k(t)B_1(\mathbf{0}).$$

We can now write the solution $X_1(t)$ of (3.1) as follows:

$$X_1(t) = U\sum_{k=0}^\infty \left(\frac{1}{k!}\Sigma_1^k(t)\right)B_1(\mathbf{0}) + U\sum_{k=0}^\infty \left(\frac{1}{k!}\Sigma_2^k(t)\right)B_1(\mathbf{0}) = U\left(e^{\Sigma_1(t)}B_1(\mathbf{0}) + e^{\Sigma_2(t)}B_1(\mathbf{0})\right)$$

It is clear that the $U$ specifies the rotation of $e^{\Sigma_1(t)}B_1(\mathbf{0}) + e^{\Sigma_2(t)}B_1(\mathbf{0})$.

**Example 3.2.** Let

$$A(t) = \begin{pmatrix} e^{-t} + 2e^{-t}(4 - \cos(t)) & e^{-t}(4 - \cos(t)) \\ e^{-t}(4 - \cos(t)) & e^{-t} \end{pmatrix}, \quad \forall t \geq t_0 = 0$$

We can see that the matrix $A(t)$ is symmetric, $\delta(t) \neq 0$, and it satisfies the Lappo-Danilevskii condition $p(t) = e^{-t}, q(t) = e^{-t}(4 - \cos(t))$, and $\gamma = 2$.

The eigenvalues of the matrix $\int_s^t A(s)ds$ are

$$\lambda_{1,2}(s,t) = \frac{1}{2}\left(\int_s^t \left(2e^{-w} + 2e^{-w}(4-\cos(w))\right) dw \pm \sqrt{8}\left|\int_s^t e^{-w}(4-\cos(w))dw\right|\right)$$

and the singular values are:

$$\sigma_1(s,t) = |\lambda_1(s,t)| = \left| \left(5 + 4\sqrt{2}\right)(e^{-s} - e^{-t}) + \left(\frac{1+\sqrt{2}}{2}\right)(e^{-s}\sin(s) - e^{-t}\sin(t)) \right.$$
$$\left. + \left(\frac{1+\sqrt{2}}{2}\right)(e^{-t}\cos(t) - e^{-s}\cos(s)) \right|,$$

and,

$$\sigma_2(s,t) = |\lambda_2(s,t)| = \left| \left(5 - 4\sqrt{2}\right)(e^{-s} - e^{-t}) + \left(\frac{1-\sqrt{2}}{2}\right)(e^{-s}\sin(s) - e^{-t}\sin(t)) \right.$$
$$\left. + \left(\frac{1-\sqrt{2}}{2}\right)(e^{-t}\cos(t) - e^{-s}\cos(s)) \right|.$$

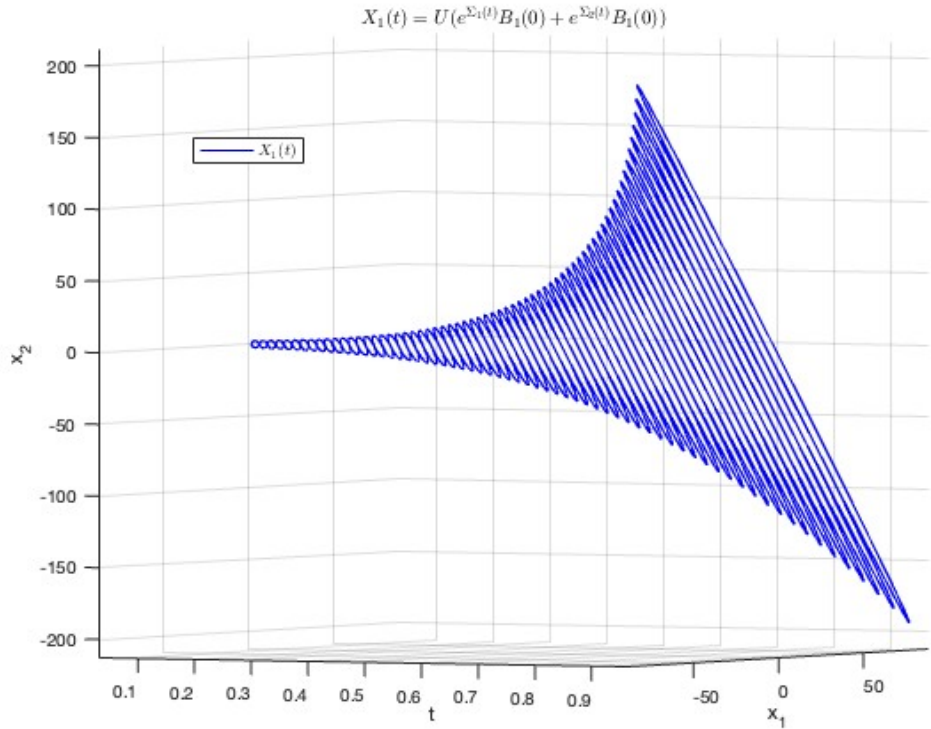

Figure 5: Solution $X_1(t), t \in [0, 1]$

Let $r(t) = e^t$, then the differential equation (3.1) possesses the solution $X_1(t) = U\left(e^{\Sigma_1(t)}B_1(\mathbf{0}) + e^{\Sigma_2(t)}B_1(\mathbf{0})\right)$ (Fig. (5)) and rotated by an angle $\theta \approx 22.5°$ specified by the matrix $U$ given by relation (3.6) as follows:

$$U = \begin{pmatrix} 0.9239 & 0.3827 \\ 0.3827 & -0.9239 \end{pmatrix}$$

**Remark 3.3.** For the initial set $X_0$, and the forcing set $F$, the following statements hold:

(i) If $X_0 = (-1)X_0$, then $X_0 \stackrel{h}{-} (-1)X_0$ exists and equals $\{\mathbf{0}\}$.

(ii) If $X_0 = (-1)X_0$ and $AX_0 = \sigma X_0$, then there exists an $\alpha > 0$ such that H-differences $X_0 \stackrel{h}{-} \alpha A X_0$ and $X_0 \stackrel{h}{-} \alpha(-1)AX_0$ both exist.

(iii) If $F \in K_c(\mathbb{R}^n)$ is such that $F = (-1)F$ and $AF = \sigma F$, then there exists an $\alpha > 0$ such that H-differences $F \stackrel{h}{-} \alpha AF$ and $F \stackrel{h}{-} \alpha(-1)AF$ exist.
Moreover, under these conditions, the solutions of differential equations with PS-derivative will also be solutions of the differential equation with BG-derivative and vice versa.

We use the facts stated in remark (3.3) in the example below, wherein we discuss in detail about the mixed solutions of an SVDE.

**Example 3.3.** Consider the initial value problem (IVP)

$$DX(t) = A(t)X(t) + B_{r(t)}(\mathbf{0}), \quad X(0) = B_1(\mathbf{0}), \quad t \in [0,1] \tag{3.7}$$

where $A(t) \in \mathbb{R}^{2\times 2}$ is given by

$$A = \begin{pmatrix} a(t)\cos(\varphi) & -a(t)\sin(\varphi) \\ a(t)\sin(\varphi) & a(t)\cos(\varphi) \end{pmatrix},$$

with $a(t)$ is continuous function for all $t \geq 0$, and $\varphi \in [0, 2\pi)$. The derivative $DX(t)$ is one of the previously considered derivatives $(D_H X(t), D_{ps} X(t), D_{bg} X(t))$ of the set-valued mapping $X(t): [0,1] \to K_c(\mathbb{R}^2)$, and $B_{r(t)}(\mathbf{0}) = \{x \in \mathbb{R}^2 : \|x\| \leq r(t)\}$ is a closed ball of radius $r(t)$ for all $t \geq 0$ centered at $\mathbf{0} \in \mathbb{R}^2$.

Since the matrix $A(t)$ can be expressed as $A(t) = a(t)R(\varphi)$, where $R(\varphi)$ is a rotation matrix, the singular values of $A(t)$ are equal for any $a(t)$ and $\varphi$, namely

$$\sigma_1(t) = \sigma_2(t) = |a(t)|.$$

We have (see [21])

$$A(t)B_1(\mathbf{0}) = a(t)R(\varphi)B_1(\mathbf{0}) = |a(t)|B_1(\mathbf{0}) = B_{|a(t)|}(\mathbf{0}).$$

By theorem (3.1), the first and the second basic solutions to the IVP (3.7) are given, respectively, by:

$$X_1(t) = \left(e^{\int_0^t |a(u)|\,du} + \int_0^t e^{\int_s^t |a(u)|\,du} r(s)\,ds\right) B_1(\mathbf{0}),$$

$$X_2(t) = \left(e^{-\int_0^t |a(u)|\,du} - \int_0^t e^{-\int_s^t |a(u)|\,du} r(s)\,ds\right) B_1(\mathbf{0}).$$

In particular, if $A(t) = tR(\varphi)$, the singular values of $A(t)$ are $\sigma_1 = \sigma_1 = |t|$, and with $r(t) = t$, the first and the second basic solutions to the IVP (3.7) are given, respectively, by:

$$X_1(t) = \left(2e^{\frac{t^2}{2}} - 1\right)B_1(\mathbf{0}), \quad \text{(Fig. (6))}, \qquad X_2(t) = \left(2e^{-\frac{t^2}{2}} - 1\right)B_1(\mathbf{0}), \quad \text{(Fig. (7))}.$$

For further discussion regarding mixed solutions, for simplicity, we consider the case where the coefficient matrix $A(t) \equiv A$ is constant. Then the matrix $A = aR(\varphi)$, and the singular values of $A$ are $\sigma_1 = \sigma_2 = |a|$, and let $r(t) = r$, also be a constant. The first and the second basic solutions in this case (see [18]), are given by:

$$X_1(t) = \left(\left(1 + \frac{r}{|a|}\right) e^{|a|t} - \frac{r}{|a|}\right) B_1(\mathbf{0}), \quad \text{(Fig. (8))},$$

$$X_2(t) = \left(\left(1 + \frac{r}{|a|}\right) e^{-|a|t} - \frac{r}{|a|}\right) B_1(\mathbf{0}), \quad t \leq \frac{1}{|a|} ln\left(1 + \frac{|a|}{r}\right), \quad \text{(Fig. (9))}.$$

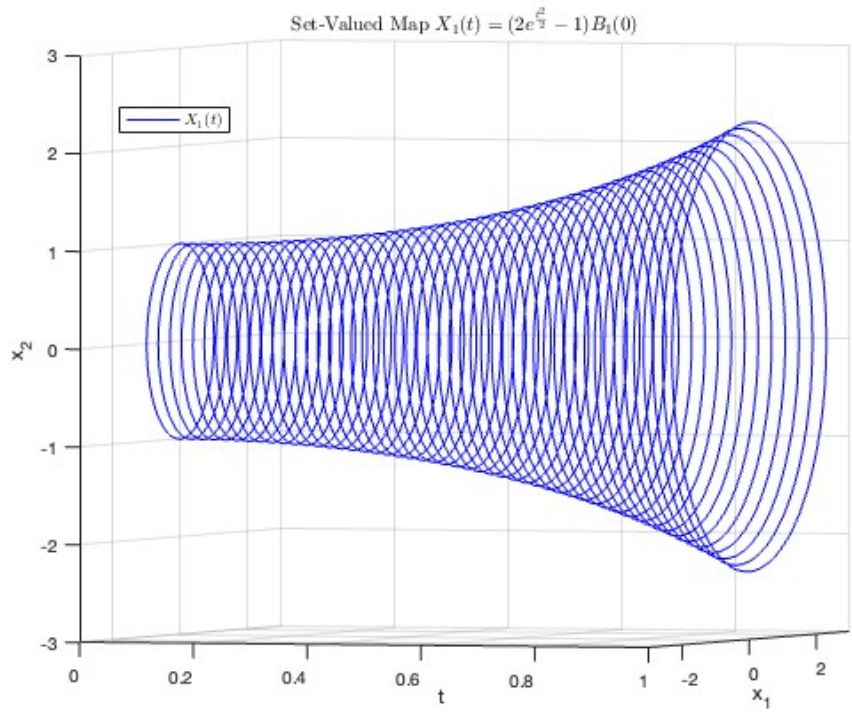

Figure 6: Solution $X_1(t), t \in [0,1]$

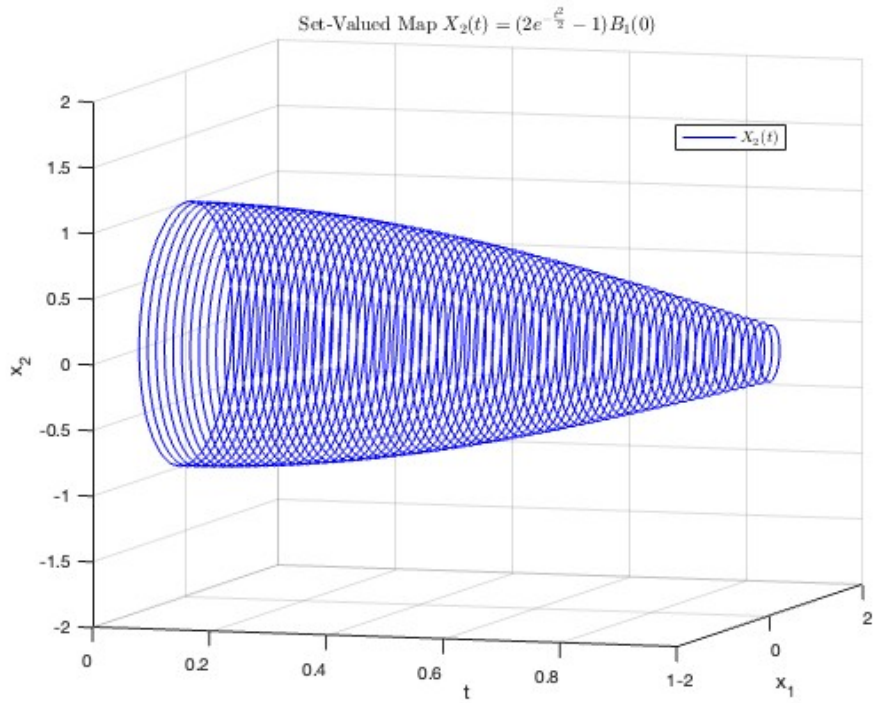

Figure 7: Solution $X_2(t), t \in [0,1]$

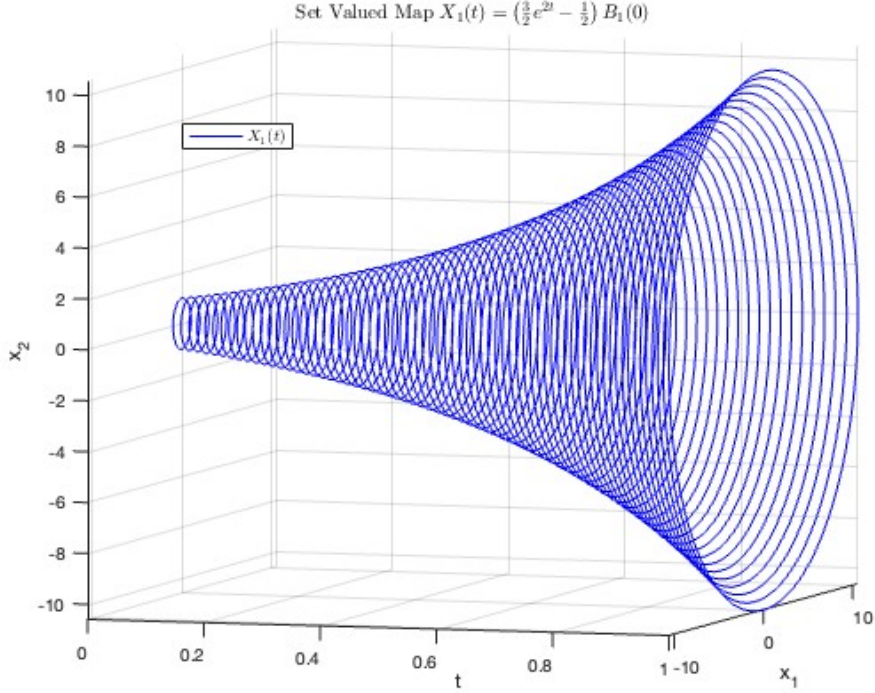

Figure 8: $a = 2, r = 1, X_1(t), t \in [0, 1]$

By [18], the second basic solution $X_2(\cdot)$ is also a solution of the corresponding integral equation:

$$X_2^{ps}(t) = X_0 \stackrel{h}{-} A \int_0^t X_2^{ps}(s) ds \stackrel{h}{-} tB_r(\mathbf{0}) \quad \text{or} \quad X_2^{bg}(t) = X_0 \stackrel{h}{-} (-1) A \int_0^t X_2^{bg}(s) ds \stackrel{h}{-} (-1) tB_r(\mathbf{0})$$

Using the basic solutions, we can now construct mixed solutions to this IVP (3.7) with PS-derivative (BG-derivative).

The set-valued mappings

$$Y_1(t) = \begin{cases} \left(\left(1 + \frac{r}{|a|}\right) e^{|a|t} - \frac{r}{|a|}\right) B_1(\mathbf{0}), & t \in [0, 0.5]; \\ \left(\left(1 + \frac{r}{|a|}\right) e^{|a|(1-t)} - \frac{r}{|a|}\right) B_1(\mathbf{0}), & t \in [0.5, 1] \end{cases}$$

and,

$$Y_2(t) = \begin{cases} \left(\left(1 + \frac{r}{|a|}\right) e^{-|a|t} - \frac{r}{|a|}\right) B_1(\mathbf{0}), & t \in [0, 0.5]; \\ \left(\left(1 + \frac{r}{|a|}\right) e^{|a|(t-1)} - \frac{r}{|a|}\right) B_1(\mathbf{0}), & t \in [0.5, 1] \end{cases}$$

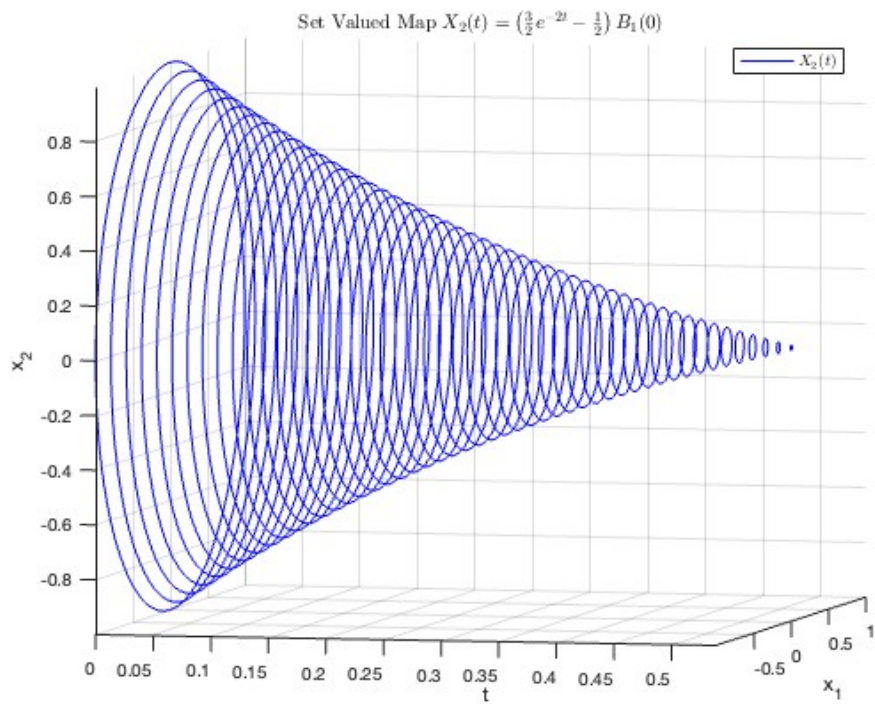

Figure 9: $a = 2,\ r = 1,\ X_2(t), t \in [0,1]$

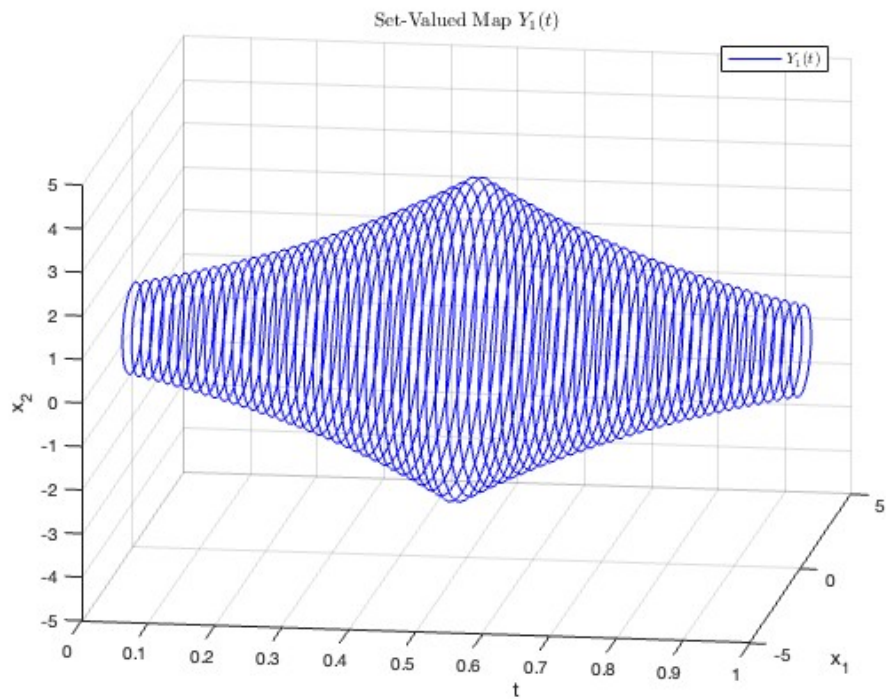

Figure 10: $a = 2,\ r = 1,\ Y_1(t), t \in [0,1]$

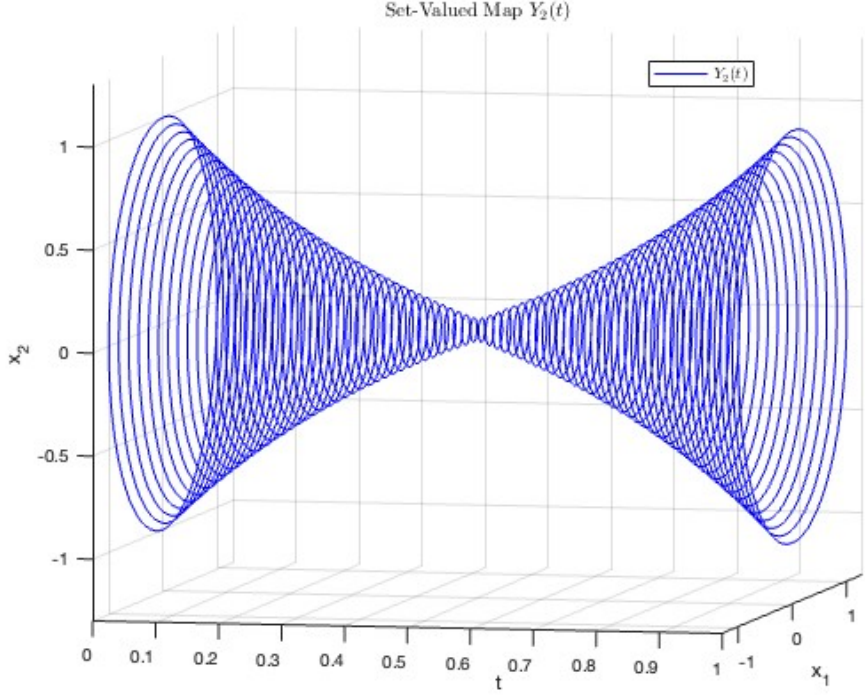

Figure 11: $a = 2$, $r = 1$, $Y_2(t), t \in [0, 1]$

are the solutions of the IVP (3.7) with the PS-derivative and also with the BG-derivative (Figs. (10) & (11)). Note that we can construct infinitely many such mixed solutions using the basic solutions. For these mixed solutions $Y(\cdot)$, the diameter function diam $(Y(\cdot))$ is not monotonically increasing or monotonically decreasing in the entire considered interval (or on the interval where they exist).

Also we note that the solution $Y_1(\cdot)$ is a solution of the integral equations:

$$Y_1^{ps}(t) = X_0 + \int_0^{u(t)} \left(AY_1^{ps}(s) + B_r(\mathbf{0})\right)ds \stackrel{h}{-} \mathcal{U}(t - 0.5) \int_{0.5}^{v(t)} \left(AY_1^{ps}(s) + B_r(\mathbf{0})\right) ds,$$

and,

$$Y_1^{bg}(t) = X_0 + \int_0^{u(t)} \left(AY_1^{bg}(s) + B_r(\mathbf{0})\right)ds \stackrel{h}{-} (-1)\mathcal{U}(t - 0.5) \int_{0.5}^{v(t)} \left(AY_1^{bg}(s) + B_r(\mathbf{0})\right) ds,$$

where $t \in [0, 1]$, $u(t) = \min\{t, 0.5\}$, $v(t) = \max\{t, 0.5\}$, $\mathcal{U}(t)$ is the Heaviside (unit) step function.

Similarly, the solution $Y_2(\cdot)$ is a solution of the integral equations:

$$Y_2^{ps}(t) = X_0 \stackrel{h}{-} \int_0^{u(t)} \left(AY_2^{ps}(s) + B_r(\mathbf{0})\right)ds + \mathcal{U}(t - 0.5) \int_{0.5}^{v(t)} \left(AY_2^{ps}(s) + B_r(\mathbf{0})\right)ds,$$

and,

$$Y_2^{bg}(t) = X_0 \stackrel{h}{-} (-1)\int_0^{u(t)} \left(AY_2^{bg}(s) + B_r(\mathbf{0})\right)ds + \mathcal{U}(t - 0.5) \int_{0.5}^{v(t)} \left(AY_2^{bg}(s) + B_r(\mathbf{0})\right)ds.$$

**Remark 3.4.** It is also worth noting that, in this example (3.3), the shape of the cross-section of the solution corresponds to the shape of the initial set, and the diameter of each section remains independent of the parameter $\varphi$.

However, if we take $A = \begin{pmatrix} a\cos(\varphi) & -b\sin(\varphi) \\ b\sin(\varphi) & a\cos(\varphi) \end{pmatrix}$, where $a, b \in \mathbb{R}$, and $\varphi \in [0, 2\pi)$, then the singular values of the matrix $A$ are equal to each other for any $a, b$ and $\varphi$ and $\sigma_1(a, b, \varphi) = \sigma_2(a, b, \varphi) = \sigma(a, b, \varphi) = \sqrt{a^2\cos^2(\varphi) + b^2\sin^2(\varphi)}$. Therefore, the matrix $A = \sigma(a, b, \varphi)R(\psi)$, where $\psi = \arccos\left(\frac{a\cos(\varphi)}{\sigma_1(a,b,\varphi)}\right)$. Furthermore, the set-valued mappings

$$X_1(t) = \left(\left(1 + \frac{r}{\sigma(a,b,\varphi)}\right)e^{\sigma(a,b,\varphi)t} - \frac{r}{\sigma(a,b,\varphi)}\right) B_1(\mathbf{0}),$$

and,

$$X_2(t) = \left(\left(1 + \frac{r}{\sigma(a,b,\varphi)}\right)e^{-\sigma(a,b,\varphi)t} - \frac{r}{\sigma(a,b,\varphi)}\right) B_1(\mathbf{0})$$

are basic solutions of the IVP (3.7) with the PS-derivative and with the BG-derivative. In this case, the size of the section of the solutions will depend on $a, b$ and $\varphi$.

## 4 Conclusion

In this paper, we studied the solutions of non-homogeneous linear SVDEs with variable matrix-valued coefficients in the space of compact convex subsets of $\mathbb{R}^2$, with the generalized derivative that includes Hukuhara, BG- and PS-derivatives. We obtained explicit analytical formulas for the basic and mixed solutions under suitable assumptions on the time-varying coefficient matrix. These results highlight that, to ensure the existence of solutions to such problems, the initial set and the forcing term must be chosen appropriately. We utilize the central symmetry and the property $(-1)B_1(\mathbf{0}) = B_1(\mathbf{0})$ of the balls centered at the origin, to obtain our results. Additional restrictions come into play while formulating such problems with arbitrary elements of $K_c(\mathbb{R}^n)$ and will require a separate study.

### Statements and Declarations

**Competing Interests:** The authors have no competing interests to declare that are relevant to the content of this article.

**Funding Information:** No funding was received during the preparation of this manuscript.

**Author Contribution:** Both authors contributed equally to the established results, read, and approved the final manuscript.

**Data Availability:** There is no data associated with this manuscript.

## References

[1] F. S. De Blasi and F. Iervolino, "Equazioni differenziali con soluzioni a valore compatto convesso," Bollettino della Unione Matematica Italiana. Series IV, vol. 2, pp. 491–501, 1969.

[2] M. Hukuhara, "Intégration des applications mesurables dont la valeur est un compact convexe," Funkcialaj Ekvacioj. Serio Internacia, vol. 10, pp. 205–223, 1967.

[3] V. Lakshmikantham, T. G. Bhaskar, and J. V. Devi, Theory of set differential equations in metric spaces. Cottenham, Cambridge: CSP, Cambridge Scientific Publ, 2006.

[4] B. Bede and S. G. Gal, "Generalizations of the differentiability of fuzzy-number-valued functions with applications to fuzzy differential equations," Fuzzy Sets and Systems, vol. 151, pp. 581–599, May 2005.

[5] L. Stefanini and B. Bede, "Generalized Hukuhara differentiability of interval-valued functions and interval differential equations," Nonlinear Analysis: Theory, Methods & Applications, vol. 71, pp. 1311–1328, Aug. 2009.

[6] Plotnikov, "Set-valued differential equations with generalized derivative," Journal of Advanced Research in Pure Mathematics, vol. 3, pp. 144–160, Jan. 2011.

[7] T. Komleva, A. Plotnikov, and N. Skripnik, "Some properties of the solutions of a linear set-valued differential equation in the space conv($\mathbb{R}^2$)," Journal of Mathematical Sciences, vol. 279, pp. 363–383, Feb. 2024.

[8] T. A. Komleva, A. V. Plotnikov, L. I. Plotnikova, and N. V. Skripnik, "Conditions for the existence of basic solutions of linear multivalued differential equations," 2021, vol. 73, May 2021.

[9] T. Komleva, L. Plotnikova, N. Skripnik, and A. Plotnikov, "Some remarks on linear set-valued differential equations," Studia Universitatis Babeș-Bolyai Mathematica, vol. 65, pp. 411–427, Sept. 2020. Number: 3.

[10] T. Komleva, A. Plotnikov, and N. Skripnik, "Existence of solutions of linear set-valued integral equations and their properties," Journal of Mathematical Sciences, vol. 277, pp. 1–13, Dec. 2023.

[11] A. Plotnikov and N. Skripnik, "Existence and Uniqueness Theorems for Generalized Set Differential Equations," International Journal of Control Science and Engineering, vol. 2, no. 1, pp. 1–6, 2012. Publisher: Scientific & Academic Publishing.

[12] A. V. Plotnikov, T. A. Komleva, and N. V. Skripnik, "Existence of basic solutions of first order linear homogeneous set-valued differential equations," Matematychni Studii, vol. 61, pp. 61–78, Mar. 2024. Number: 1.

[13] N. A. Gasilov, "Solving a system of linear differential equations with interval coefficients," Discrete & Continuous Dynamical Systems - B, vol. 26, no. 5, p. 2739, 2021.

[14] E. Amrahov, A. Khastan, N. Gasilov, and A. G. Fatullayev, "Relationship between bede–gal differentiable set-valued functions and their associated support functions," Fuzzy Sets and Systems, vol. 295, pp. 57–71, July 2016.

[15] I. Atamas and V. Slyn'ko, "Liouville formula for some classes of differential equations with the hukuhara derivative," Differential Equations, vol. 55, pp. 1407–1419, Nov. 2019.

[16] I. V. Atamas' and V. I. Slyn'ko, "Stability of Fixed Points for a Class of Quasilinear Cascades in the Space conv n," Ukrainian Mathematical Journal, vol. 69, pp. 1354–1369, Feb. 2018. Company: Springer Distributor: Springer Institution: Springer Label: Springer Number: 9 Publisher: Springer US.

[17] I. Atamas, S. Dashkovskiy, and V. Slynko, "Impulsive Input-to-State Stabilization of an Ensemble," Set-Valued and Variational Analysis, vol. 31, p. 25, Aug. 2023.

[18] U. M. R. Epuganti and G. B. Tenali, "Solutions of non-homogeneous linear set-valued differential equations," Nonlinear Analysis: Real World Applications, vol. 87, p. 104411, Feb. 2026.

[19] N. A. Perestyuk, V. A. Plotnikov, A. M. Samoilenko, and N. V. Skripnik, "Differential equations with impulse effects: Multivalued right-hand sides with discontinuities," in Differential Equations with Impulse Effects, De Gruyter, July 2011.

[20] V. G. Boltyanski and J. Jerónimo Castro, "Centrally symmetric convex sets," Journal of Convex Analysis, vol. 14, no. 2, pp. 345–351, 2007.

[21] G. E. Forsythe and C. B. Moler, Computer Solution of Linear Algebraic Systems. Prentice-Hall, Inc., first edition ed., Jan. 1967.

[22] A. V. Plotnikov, "Differentiation of Multivalued Mappings. T-Derivative," Ukrainian Mathematical Journal, vol. 52, pp. 1282–1291, Aug. 2000.

[23] H. T. Banks and M. Q. Jacobs, "A differential calculus for multifunctions," Journal of Mathematical Analysis and Applications, vol. 29, pp. 246–272, Feb. 1970.

[24] Y. Chalco-Cano, H. Román-Flores, and M. D. Jiménez-Gamero, "Generalized derivative and $\pi$-derivative for set-valued functions," Information Sciences, vol. 181, pp. 2177–2188, June 2011.

[25] N. Erougin, "Reducible systems," Trav. Inst. Math. Stekloff, vol. 13, p. 95, 1946.